\documentclass[12pt]{amsart}
\usepackage{amsfonts, amsmath, amssymb, amscd, latexsym}

\newif\ifpdf
    \ifx\pdfoutput\undefined
    \pdffalse % we are not running PDFLaTeX
    \else
    \pdfoutput=1 % we are running PDFLaTeX
    \pdftrue
    \fi

    \ifpdf
    \usepackage[pdftex]{graphicx}
    \else
    \usepackage{graphicx}
    \fi

\newtheorem{thm}{Theorem}[section]
\newtheorem{prop}[thm]{Proposition}

\newtheorem{cor}[thm]{Corollary}

\newtheorem{defn}[thm]{Definition}

\newtheorem{lemma}[thm]{Lemma}
\newtheorem{corollary}[thm]{Corollary}

\begin{document}

\title{Combinatorial and metric properties of Thompson's group $T$}

\author{Jos\'e Burillo}
\address{
Departament de Matem\'atica Aplicada IV, Universitat Polit\'ecnica de Catalunya, Escola Polit{\`e}cnica
Superior de Castelldefels, 08860 Castelldefels, Barcelona, Spain} \email{burillo@mat.upc.es}

\author{Sean Cleary}
\address{Department of Mathematics,
The City College of New York \& The CUNY Graduate Center, New York, NY 10031}
\email{cleary@sci.ccny.cuny.edu}

\author{Melanie Stein}
\address{Department of Mathematics, Trinity College, Hartford, CT
06106} \email{melanie.stein@trincoll.edu}

\author{Jennifer Taback}
\address{Department of Mathematics, Bowdoin College, Brunswick, ME
04011} \email{jtaback@bowdoin.edu}

\thanks{The first, second and fourth authors acknowledge support from
NSF International Collaboration grant DMS-0305545 and are grateful for the hospitality of the Centre de
Recerca Matem\`atica.}

%\subjclass[2000]{Primary ; Secondary }
%\thanks{}
%\keywords{}

\date{\today}

%\dedicatory{}

%%% ----------------------------------------------------------------------

\begin{abstract}
We discuss metric and combinatorial properties of Thompson's group $T$, including  normal forms for
elements and unique tree pair diagram representatives.  We relate these properties to those of Thompson's group
$F$ when possible, and highlight combinatorial differences between the two groups.  We define a set of
unique normal forms for elements of $T$ arising from minimal factorizations of elements into natural pieces.  We show that the number of carets in a reduced representative of  an element  of $T$ estimates the word length, and
that  $F$ is undistorted in $T$. We describe how to recognize torsion
elements in $T$.
\end{abstract}

\maketitle

\section{Introduction}

Thompson's groups $F$, $T$ and $V$ are a remarkable family of infinite, finitely-presentable groups
studied for their own properties as well as for their connections with questions in logic,
homotopy theory, geometric group theory and the amenability of discrete groups.

Cannon, Floyd and Parry give an excellent introduction to these groups in \cite{cfp}.  These three groups
can be viewed either algebraically, combinatorially, or analytically.  Algebraically, each has both
finite and infinite presentations. Geometrically, an element in each group can be viewed as a {\em tree
pair diagram}; that is, as a pair of finite binary rooted trees with the same number of leaves, with a
numbering system pairing the leaves in the two trees.   Analytically, an element of each group can be
viewed as a piecewise-linear self map of the unit interval:
\begin{itemize}
\item in $F$ as a piecewise linear homeomorphism,
\item in $T$ as a homeomorphism of the unit interval with the endpoints identified, and thus of $S^1$,
\item in $V$ as a right-continuous bijection which is locally orientation preserving.
\end{itemize}

%
%
%we can regard the rooted binary trees as giving
%dyadic subdivisions of the unit interval, and the piecewise-linear interpolation of these
%subdivisions yields a map from the unit interval to itself.   If $w \in F$, then the leaves in
%both trees in a diagram representing $w$ are numbered from left to right and this pairing gives
%rise to an orientation-preserving homeomorphism of the interval $[0,1]$.
% If $w \in T$, then the leaves in the source tree are numbered from
%left to right, and in the target tree any leaf can be numbered zero, with the rest of the leaves
%labelled cyclically in increasing order.  This will give rise to an orientation-preserving
%homeomorphism of the unit circle $S^1$, regarded as  $[0,1]$ with the endpoints identified.  If $w
%\in V$, again we number the leaves in the source tree from left to right, but the leaves in the
%target tree are numbered in any order, and paired according to this numbering.  The resulting
%self-map of $[0,1]$ in this case is not necessarily a homeomorphism, but rather a
%right-continuous bijection from the interval to itself, which is locally orientation-preserving.
%

Thompson's group $F$ in particular has been studied extensively. The group $F$ has a standard infinite
presentation in which every element has a unique normal form, and a standard two-generator finite
presentation.  Fordham \cite{blakegd} presented a method of computing the word length of $w \in F$ with
respect to the standard finite generating set directly from a tree pair diagram representing $w$.
Regarding $F$ as a diagram group, Guba \cite{gubagrowth} also obtained an effective geometric method for
computing the word metric with respect to the standard finite generating set.  Belk and Brown
\cite{belkbrown} have similar results which arise from viewing elements of $F$ as forest diagrams. 

%The
%amenability of $F$ is one open question which has motivated continuing research into $F$ for many years.
%Either answer to this question would be interesting, for if $F$ is not amenable it would be  a familiar
%finitely presented non-amenable group without free non-abelian subgroups, whereas if $F$ is amenable it
%would be a familiar finitely presented, amenable, but not elementary-amenable group.

In this paper, we discuss analogues for $T$ of some properties of $F$, using all
three of the descriptions of $T$:  algebriac, geometric and analytic. We  begin by
desribing unique normal forms for elements which arise from their  reduced tree pair descriptions.
We consider metrically how $F$ is contained as a subgroup of $T$, and show that the number of carets in a
reduced tree pair diagram representing $w \in T$ estimates the word length of $w$ with respect to a
particular generating set.  Thus $F$ is quasi-isometrically embedded in $T$. Furthermore,
we show that there are families of words in $F$ which are isometrically embedded in $T$ with respect to an alternate finite generating set.  The groups $T$ and $V$, unlike $F$, contain
torsion elements, and we describe how to recognize these torsion elements from their tree pair diagrams. Finally, we show that every  torsion element of $T$ is  conjugate to a power to a generators of $T$
and that the subgroup of rotations in $T$ is quasi-isometrically embedded. 
\section{Background on Thompson's groups $F$ and $T$}

\subsection{Presentations and tree pair diagrams}

Thompson's groups $F$ and $T$ both have representations as groups of piecewise-linear homeomorphisms. The
group $F$ is the group of orientation-preserving homeomorphisms of the interval $[0,1]$, where each
homeomorphism is required to have only finitely many discontinuities of slope, called {\em breakpoints},
have slopes which are powers of two and have the coordinates of the breakpoints all lie in the set of dyadic
rationals. Similarly, the group $T$ consists of orientation-preserving homeomorphisms of the circle $S^1$
satisfying the same conditions where we represent the circle $S^1$ as the unit interval $[0,1]$ with the
two endpoints identified.

Cannon, Floyd and Parry give an excellent introduction to Thompson's groups $F$,  $T$ and $V$ in
\cite{cfp}. We refer the reader to this paper for full details on results mentioned in this section.
Since more readers have some familiarity with $F$ than with $T$, we first give a very brief review of the
group $F$, and then a slightly more detailed review of $T$. Algebraically, $F$ has well known infinite
and finite presentations. With respect to the infinite presentation
$$
\langle x_i, i\geq 0\, |\,x_jx_i=x_ix_{j+1}, i<j\rangle
$$
for $F$, group elements have simple normal forms which are unique. It is easy to see that $F$ can be
generated by $x_0$ and $x_1$, which form the standard finite generating set for $F$, and yield the finite
presentation

$$
\langle x_0,x_1\,|\,[x_0x_1^{-1},x_0^{-1}x_1x_0],[x_0x_1^{-1},x_0^{-2}x_1x_0^2]\rangle.
$$

A geometric representation for an element $w$ in $F$ is a tree pair diagram, as discussed in \cite{cfp}.
A {\em tree pair diagram} is  a pair of finite rooted binary trees with the same number of leaves. By
convention, the leaves of each tree are thought of as being numbered from $0$ to $n$ reading from left to
right. A node of the tree together with its two downward directed edges is called a {\em caret}. The {\em
left side} of the tree consists of the root caret, and all carets connected to the root by a path of left
edges;  the {\em right side} of the tree is defined analogously.  A caret is called a {\em left caret} if
its left leaf lies on the left side of the tree.  A caret is called a {\em right caret} if it is not the
root caret and its right leaf lies on the right side of the tree.  All other carets are called {\em
interior}. A caret is called {\em exposed} if it contains two leaves of the tree. For $w \in F$, we write
$w = (T_-,T_+)$ to express $w$ as a tree pair diagram, and refer to $T_-$ as the {\em source} tree and
$T_+$ as the {\em target} tree. These trees arise naturally from the interpretation of $F$ as a group of
homeomorphisms. Thinking of $w$ as a homeomorphism of the unit interval, the source tree represents a
subdivision of the domain into subintervals of width $1/2^n$ for varying values of $n$, and the target tree represents another such
a subdivision of the range. The homeomorphism then maps the $i^{th}$ subinterval in the domain linearly
to the $i^{th}$ subinterval in the range.

A tree pair diagram representing $w$ in $F$ is not unique. A new diagram can always be produced from a given tree pair
diagram representing $w$ simply by adding carets to the $i$th leaf of both trees. We impose a natural
reduction condition: if $w = (T_-,T_+)$ and both trees contain a caret with two exposed leaves numbered
$i$ and $i+1$, then we remove these carets, thus forming a representative for $w$ with fewer carets and
leaves.
 A tree pair diagram which admits no such reductions is called a {\em reduced tree pair diagram},
and any element of $F$ is represented by a unique reduced tree pair diagram. When we write $w =
(T_-,T_+)$ below, we are assuming that the tree pair diagram is reduced unless otherwise specified.

The group $T$ also has both a finite and an infinite presentation. The infinite presentation is given by
two families of generators, $\{x_i,i\ge 0\}$, the same generators as in the infinite presentation of $F$,
a family $\{c_i,i\ge0\}$ of torsion elements, and the following relators:
\begin{enumerate}
\item $x_jx_i=x_ix_{j+1}$, if $i<j$ \item $x_kc_{n+1}=c_{n}x_{k+1}$,
if $k<n$ \item $c_nx_0=c_{n+1}^2$ \item $c_n=x_nc_{n+1}$ \item$c_n^{n+2}=1$.
\end{enumerate}

\begin{figure}[h]
\includegraphics[width=5in]{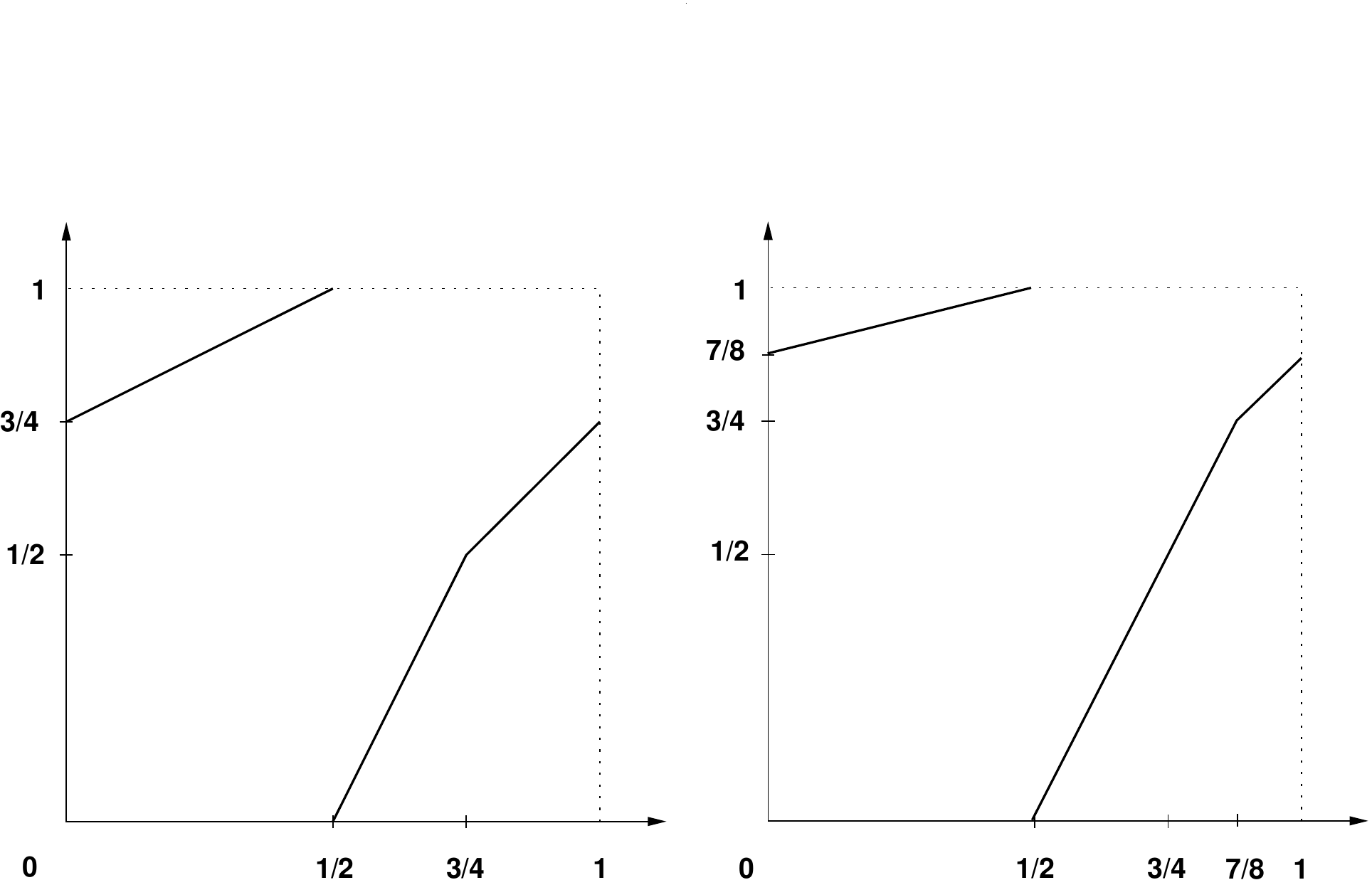}\\
\caption{ The graphs of the homeomorphisms corresponding to the elements $c_1$ and
$c_2$.\label{fig:c1c2}}
\end{figure}

This new family of generators $c_n$ (of order $n+2$), is simple to describe. The generator $c_n$ corresponds to the
homeomorphism of the circle obtained as follows. Both domain and range can be thought of as the unit
interval with the endpoints identified. We subdivide the interval into $n+1$ subintervals by successively
halving the rightmost subinterval; or in other words inserting endpoints at   $\frac{1}{2}, \frac{3}{4},
\ldots ,\frac{2^{n+1}-1}{2^{n+1}}$. Then the homeomorphism maps $[0,1/2]$ linearly to
$[\frac{2^{n}-1}{2^{n}},\frac{2^{n+1}-1}{2^{n+1}}]$, and so on around each circle. For example, the
element $c$ corresponds to the homeomorphism of $S^1$ given by
$$
c(t)=\left\{\begin{array} {ll}
\frac12t+\frac34&\text{if }0\le t<\frac12\\
2t-1&\text{if }\frac12\le t<\frac34\\
t-\frac14&\text{if }\frac34\le t\le1
\end{array}\right.
$$
Figure \ref{fig:c1c2} shows the graphs of the homeomorphisms corresponding to $c_1$ and $c_2$.

Using the first three  relators, we see that only the generators $x_0$, $x_1$ and $c_1$ are required to
generate the group, since the other generators can be obtained from these three.  In the following,  we
will use $c$ to denote the generator $c_1$.  The group $T$ is finitely presented using the following
relators,with respect to the  finite generating set $\{x_0,x_1,c\}$:

\begin{enumerate}
\item $[x_0x_1^{-1},x_0^{-1}x_1x_0]=1$ \item
$[x_0x_1^{-1},x_0^{-2}x_1x_0^2]=1$ \item $x_1c_3=c_2x_2$, (that is
$x_1(x_0^{-2}cx_1^{-2})=(x_0^{-1}cx_1)(x_0^{-1}x_1x_0)$)  \item $c_1x_0=c_2^2$, (that is,
$cx_0=(x_1^{-1}cx_0)^2$)  \item $x_1c_2=c$, (that is, $x_1(x_0^{-1}cx_1)=c) $ \item $c^3=1$.
\end{enumerate}

As with Thompson's group $F$, we will frequently work with the more convenient infinite set of generators
when constructing normal forms for elements and performing computations in the group.  We will need to
express elements with respect to a finite generating set when discussing word length.
There are two natural finite generating sets for $T$, both extending the standard
finite generating set for $F$.  The first and the one that we use primarily below is
the generating set $\{x_0, x_1, c_1\}$ used in the finite presentation above.  In Section \ref{isomembed} for the purposes of counting carets carefully, we also use the generating set $\{x_0, x_1, c_0\}$, which has the
advantage that the tree pair diagram for $c_0$ has only one caret, as opposed to $c_1$,
which has two carets, at the expense of slightly more complicated relators.

Just as for $F$, tree pair diagrams serve as efficient representations for elements of $T$. However,
since elements of $T$ represent homeomorphisms of the circle rather than the interval, the tree pair
diagram must also include a bijection between the leaves of the source tree and the leaves of the target
tree to fully encode the homeomorphism. Since this bijection can at most cyclically shift the leaves, it
is determined by the image of the leftmost leaf in the source tree. Since by convention this leaf in the
source tree is already thought of as leaf $0$, this information is recorded by writing a $0$ under the
image leaf in the target tree. Hence, for $w \in T$, a {\em marked tree pair diagram} representing $w$ is
a pair of finite rooted binary trees with the same number of leaves, together with a mark (the numeral 0)
on one leaf of the second  tree. As usual, we write $w = (T_-,T_+)$ to express $w$ as a tree pair
diagram, and refer to $T_-$ as the {\em source} tree and $T_+$ (the one with the mark) as the {\em
target} tree.  We remark that to extend this to $V$, since now the bijection of the subintervals may
permute the order in any way, the marking required on the target tree to record the bijection consists of
a number on every leaf of the target tree. Just as for $F$, there are many possible tree pair
diagrams for each element of $T$, which can be obtained by adding carets to the corresponding leaves in the source and target trees in the diagrams. However, when adding the carets, placement is guided by
the marking. The leaves of the source tree are thought of as numbered from $0$ to $n$ reading from left
to right, whereas the marking of the target tree specifies where leaf number $0$ of that tree is, and
other leaves are numbered from $1$ to $n$ reading from left to right cyclically wrapping back to the left
once you reach the rightmost leaf. With this numbering in mind, carets can be added as before to leaf $i$
of both trees.  If $i \neq 0$, the mark stays where it is. Otherwise, if $i=0$, the mark on the new
target tree is placed on the left leaf of the added caret. So for $T$, we have a similar reduction
condition: if $w = (T_-,T_+)$ and both trees contain a caret with two exposed leaves numbered $i$ and
$i+1$, then we remove these carets and renumber the leaves, moving the mark if needed, thus forming a representative for $w$ with
fewer carets and leaves.

\begin{figure}
\includegraphics[width=5in]{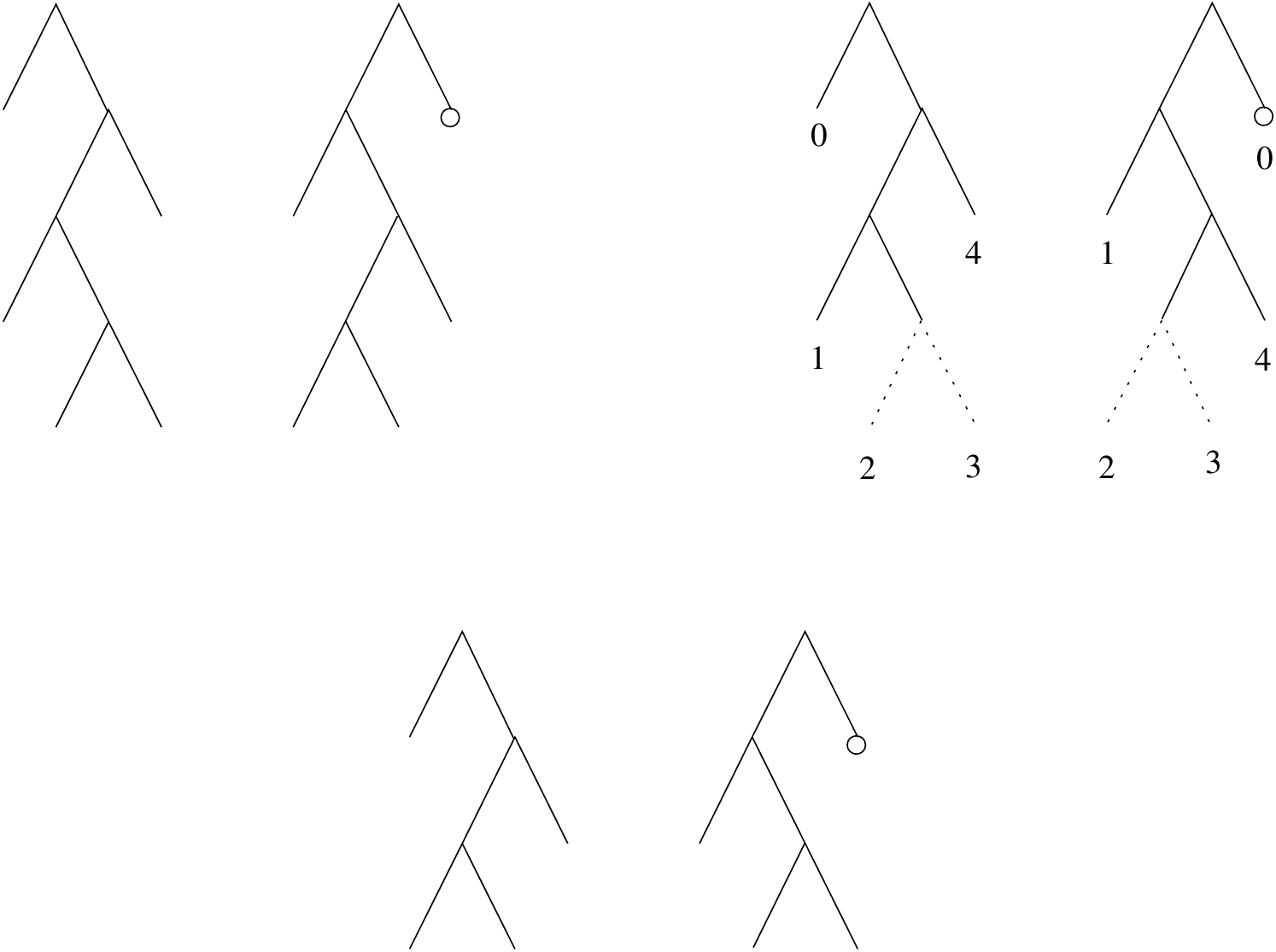}\\
\caption{An example of a caret reduction in a tree pair diagram representing an element of $T$. The diagram on the top left is reducible;  The two dotted carets on the top right are paired  with each other, since the numbering is identitcal. 
The resulting reduced diagram is shown on the bottom.\label{fig:reduction}}
\end{figure}

 A tree pair diagram which admits no such reductions is again called a {\em reduced tree pair diagram},
and any element of $T$ is represented by a unique reduced tree pair diagram. In $T$ as well as in $F$,
when we write $w = (T_-,T_+)$ below, we are assuming that the tree pair diagram is reduced unless
otherwise specified. Checking whether or not a tree pair diagram is reduced is slightly more difficult in
$T$ than in $F$. The process of checking for possible reductions is illustrated in Figure
\ref{fig:reduction}. A marked tree pair diagram for an element of $T$ is shown on the top left of Figure 
\ref{fig:reduction}. In the
top right tree pair diagram of Figure \ref{fig:reduction}, the underlying numbering of the leaves of both trees determined by the
marking is written explicitly, revealing the reducible carets. The bottom tree pair diagram shows the
resulting reduced diagram.

Note that the torsion generators $c_i$ have particularly simple tree pair diagrams. In the diagram for
$c_i$, both source and target trees consist of the root caret plus $i$ right carets. The mark $0$ is
placed on the rightmost leaf of the target tree.

\begin{figure}
\includegraphics[width=5.2in]{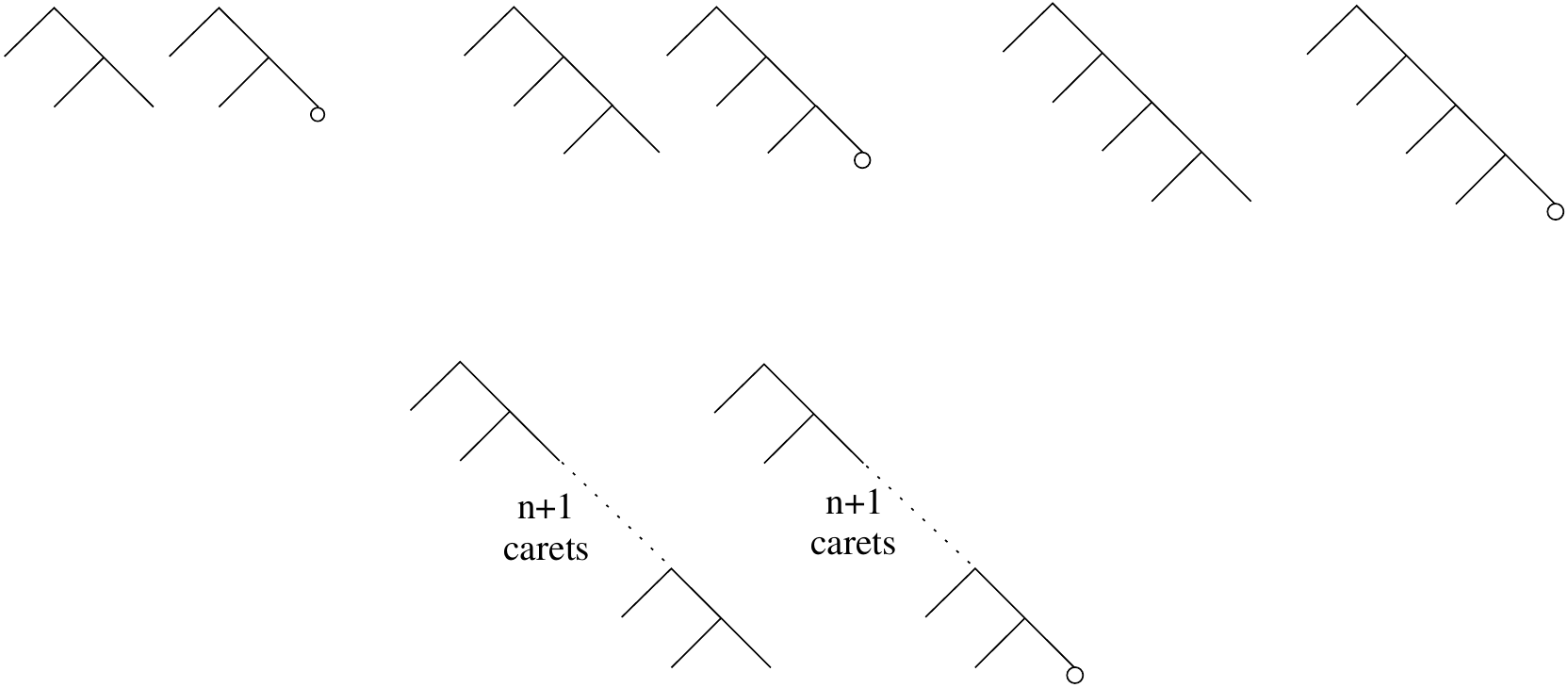}\\
\caption{Tree pair diagrams representing the elements $c_1$, $c_2$ and $c_3$ on top, plus $c_n$ on the
bottom.\label{fig:cs}}
\end{figure}

Figure \ref{fig:cs} shows the the tree pair diagrams of the first three generators $c_1, c_2$, and $c_3$, together with a
general $c_n$.  The generator $c_0$ is merely a pair of single caret trees, with the mark on the rightmost leaf of the target tree.

Whether $w \in F$ or $w \in T$, we denote the number of carets in either tree of a tree pair diagram
representing $w$ by $N(w)$.  When $p$ is a word in the generators of $F$ or $T$, then $p$ represents an
element $w$ in either $F$ or $T$, and we write $N(p)$ interchangeably with $N(w)$.

\subsection{Group Multiplication in $F$ and $T$}\label{multiplication}
Group multiplication in $F$ and $T$ corresponds to composition of
homeomorphisms, which we can interpret on the level of tree pair
diagrams as well.  First, we consider $u,v \in F$, where $u =
(T_-,T_+)$ and $v = (S_-,S_+)$.  To compute the tree pair diagram
corresponding to the product $vu$, we create unreduced
representatives $(T'_-,T'_+)$ and $(S'_-,S'_+)$ of the two elements
in which $T'_+ = S'_-$.  Then the product is represented by the
possibly unreduced tree pair diagram $(T'_-,S'_+)$. The
multiplication is written following the conventions on composition
of homeomorphisms, so the product $vu$ has as a source diagram that
of $u$, and as a target diagram that of $v$.  That is,  the diagram on
the left is the source of $u$ and the diagram on the right is the
target of $v$.

\begin{figure}
\includegraphics[width=5.2in]{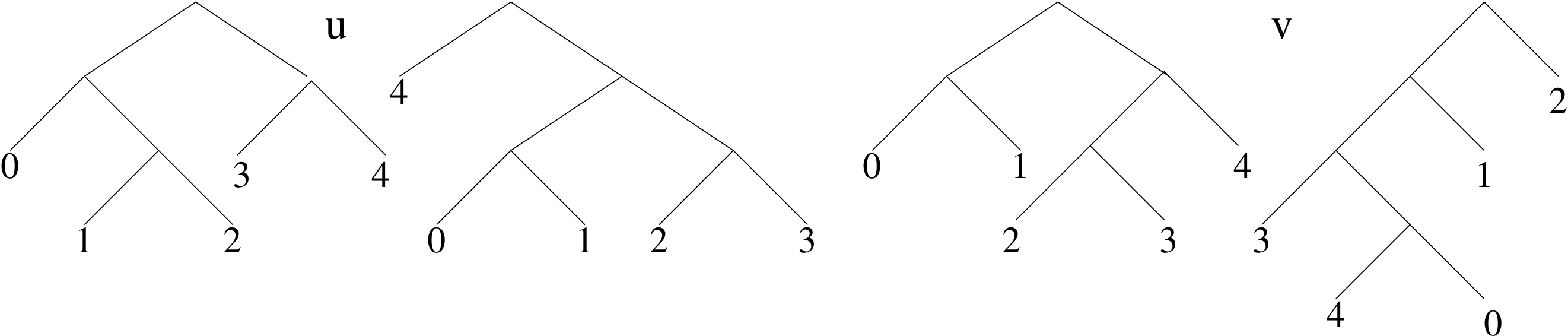}\\
\caption{ The tree pair diagram for sample elements  $u$  and $v$ in $T$.\label{fig:mult1}}
\end{figure}

\begin{figure}
\includegraphics[width=5.2in]{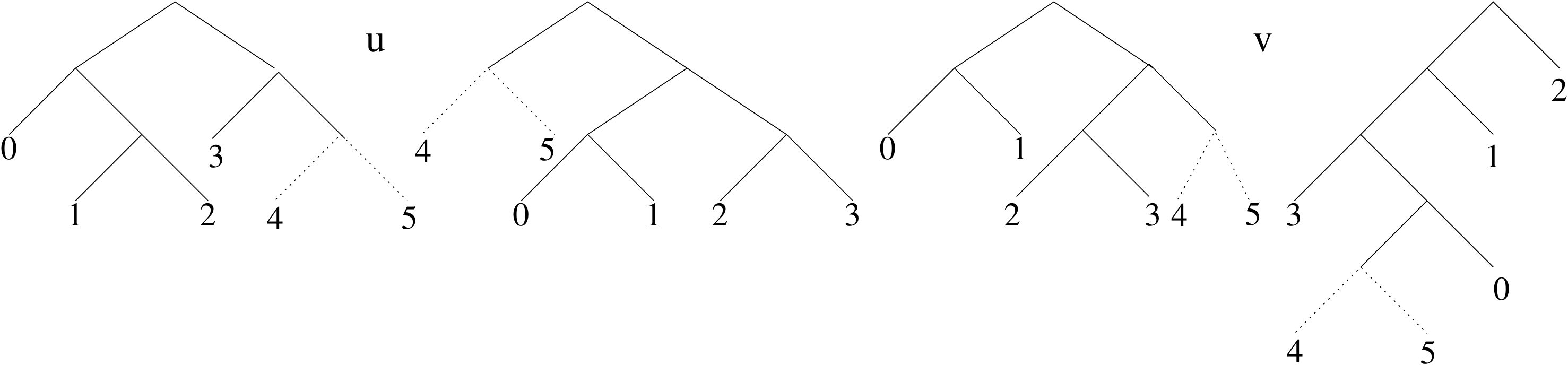}\\
\caption{ Unreduced versions of  $u$  and $v$ necessary for the multiplication $vu$ in $T$, with carets added to
perform the multiplication indicated with dashes. Now the target tree of $u$ has the same shape as the
source tree of $v$, allowing the composition.\label{fig:mult2}}
\end{figure}

\begin{figure}
\includegraphics[width=3in]{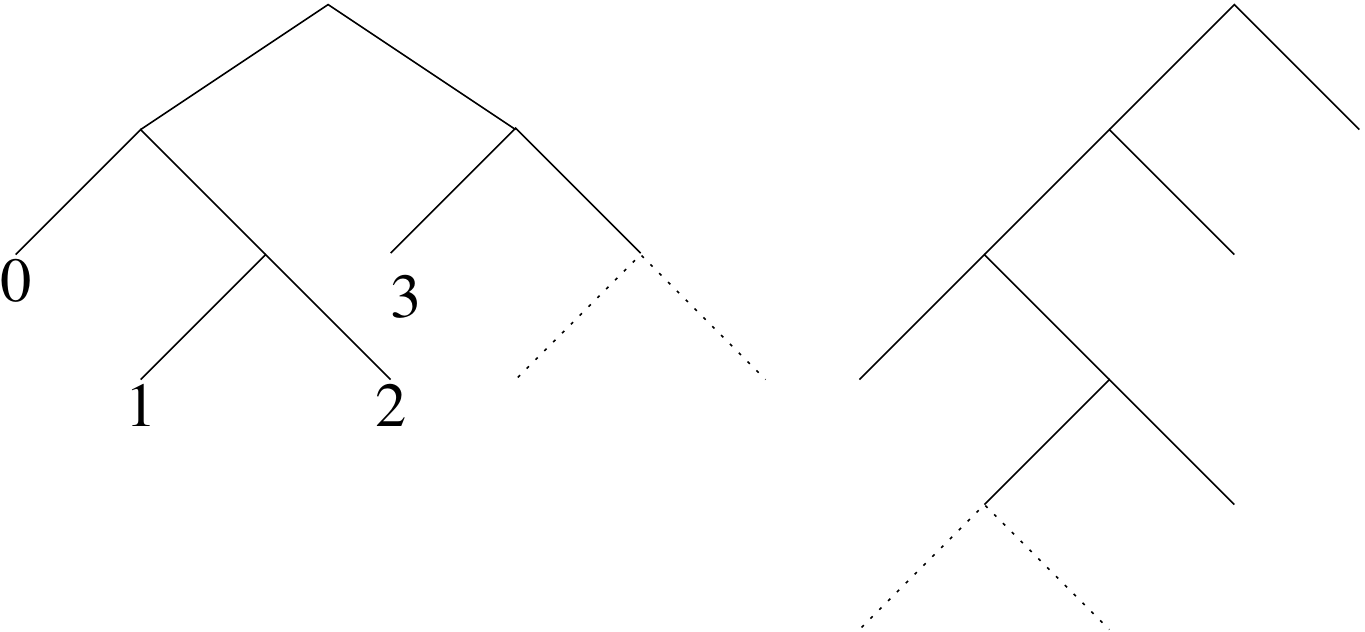}\\
\caption{ The tree pair diagram representing the product $vu$
obtained from Figure \ref{fig:mult2}. The dotted carets must be
erased to find the reduced diagram.\label{fig:mult3}}
\end{figure}

To multiply tree pair diagrams representing elements of $T$ we follow a similar procedure.  We let $u,v
\in T$, where $u = (T_-,T_+)$ and $v = (S_-,S_+)$.  To compute the tree pair diagram corresponding to the
product $vu$, we create unreduced representatives $(T'_-,T'_+)$ and $(S'_-,S'_+)$ of the two elements in
which $T'_+ = S'_-$ as trees. The product $vu$ will be represented by the pair $(T'_-,S'_+)$ of trees. To
decide which leaf in $S'_+$ to mark with the zero, we just note that it should be the leaf which is
paired with the zero leaf in $T'_-$. To identify this leaf, we find the zero leaf in $T'_+$. Since
$T'_+=S'_-$ as trees, this leaf viewed as a leaf in $S'_-$ will be labelled $m$. Then the leaf labelled
$m$ in $S'_+$ will be the new zero leaf in the tree pair diagram $(T'_-,S'_+)$ for $vu$. Alternately, we
can follow the composition in both pairs of trees to see how the leaves are paired. This newly
constructed tree pair diagram will represent $vu$ and is not necessarily reduced. For an example of this
multiplication, see Figures \ref{fig:mult1}, \ref{fig:mult2} and \ref{fig:mult3}.

%
%Finally, an interesting relation which we will use here, is the
%one between trees and elements of $F$ in particular form, as is
%established in \cite{CFP} and \cite{Blake}. Given a positive
%element
%$$
%x_{i_1}^{r_1}x_{i_2}^{r_2}\ldots x_{i_n}^{r_n}
%$$
%in $F$, where $i_1<i_2<\ldots<i_n$, corresponds to a given tree.
%The element $x_{i_k}^{r_k}$ corresponds to a sequence of $r_k$
%carets, each a left child of its predecessor, which does not reach
%the right side of the tree, and where the left leaf of the last
%caret is labelled $i_k$. This correspondence is best understood
%with an example.
%%
%% Example goes here.
%%

\section{Words and diagrams}

\subsection{Normal forms and tree pair diagrams in $F$}

With respect to the infinite presentation for $F$ given above, every element of $F$ has a unique normal
form.  Any $w$ in $F$ can be written in the form
$$w=x_{i_1}^{r_1} x_{i_2}^{r_2}\ldots x_{i_k}^{r_k}
x_{j_l}^{-s_l} \ldots x_{j_2}^{-s_2} x_{j_1}^{-s_1} $$ where $r_i, s_i >0$, $0 \leq i_1<i_2 \ldots < i_k$
and $0 \leq j_1<j_2 \ldots < j_l$. However, this expression is not unique. Uniqueness is guaranteed by
the addition of the following condition: when
 both $x_i$ and $x_i^{-1}$
occur in the expression, so does at least one of $x_{i+1}$ or  $x_{i+1}^{-1}$, as discussed by Brown and Geoghegan
\cite{bg:thomp}.   When we refer to elements of $F$ in normal form, we mean this unique normal form.

If the normal form for $w \in F$ contains no generators with negative exponents, we refer to $w$ as a
{\em positive word} and similarly, we say a normal form represents a {\em negative word} if there are no
generators with positive exponents.

We call any word which has the form $$w=x_{i_1}^{r_1} x_{i_2}^{r_2}\ldots x_{i_k}^{r_k} x_{j_l}^{-s_l}
\ldots x_{j_2}^{-s_2} x_{j_1}^{-s_1} $$ where $r_i, s_i >0$, $0 \leq i_1<i_2 \ldots < i_k$ and $0 \leq
j_1<j_2 \ldots < j_l$, a word in {\em pq form}, where $p$ is the positive part of the normal form and $q$
the negative part. The normal form for an element of $F$ is the shortest word among all words in $pq$
form representing the given element.

\begin{figure}
\includegraphics[width=3in]{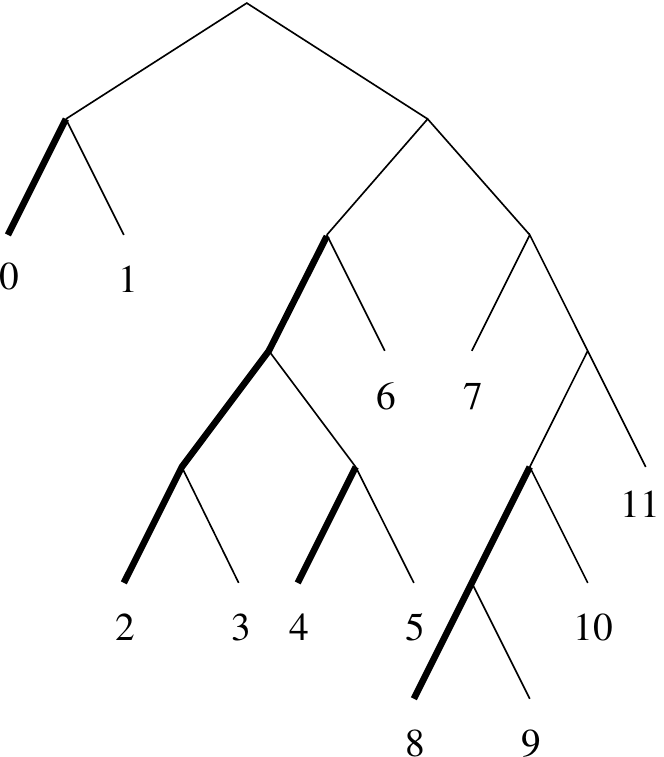}\\
\caption{Computing leaf exponents. The thick edges indicate edges
which contribute to non-zero leaf exponents. If a leaf labelled $i$ has
$r_i$ thick edges (a path of $r_i$ left edges
going up without reaching the right side of the tree) then the $i$-th leaf
exponent is $r_i$ and the generator appearing in the normal form is
$x_i^{r_i}$. This single tree $T$ pictured above is the target
tree of the tree pair diagram $(R,T)$, where $R$ is  the all-right tree with 12 leaves,
and has leaf exponents 1,0,3,0,1,0,0,0,2,0,0, and 0 for the leaves 0-11 in order.
The tree pair diagram $(R,T)$ represents the element  $x_0x_2^3x_4x_8^2$.
\label{fig:leafexp}}
\end{figure}

To any (not necessarily reduced) tree pair diagram $(T_-,T_+)$ for an element of $F$ we may associate a
word in $pq$ form representing the element, using the {\em leaf exponents} in the target and source
trees.  When the leaves of a  finite rooted binary tree are numbered from left to right, beginning with
zero, the leaf exponent of leaf $k$ is the integer length of the longest path consisting only
 of left edges of carets
which originates at leaf $k$ and does not reach the right side of the tree.  A tree pair diagram then
gives the word
$$
x_{i_1}^{r_1}x_{i_2}^{r_2}\ldots x_{i_n}^{r_n}x_{j_m}^{-s_m}\ldots x_{j_2}^{-s_2}x_{j_1}^{-s_1}
$$
precisely when leaf $i_k$ in $T_+$ has exponent $r_k$, leaf $j_k$ in
$T_-$ has leaf exponent $s_k$, and generators which do not appear in
the word correspond to leaves with exponent zero. We think of this
word as the $pq$ factorization of the element given by the
particular tree pair diagram. We call a tree an {\em all-right tree}
if it consists of a root caret together with only right carets. Note
that if we let $R$ be the all-right tree with the same number of
carets as $T_-$ or $T_+$, then $(T_-,R)$ is a diagram for the word $q$ and
$(R,T_+)$ is a diagram for word $p$. On the other hand, any word in $pq$
form can be translated into a tree pair diagram. It can be obtained
by taking diagrams for $p$ (respectively $q$), which will have all
right source (respectively target) trees. Then, if one diagram has
fewer carets, one adds right carets to its all-right tree, and of
 a corresponding path of right carets to its other tree, to make both
diagrams have exactly the same all right tree. Furthermore, under
this correspondence for $F$, reduced tree pair diagrams correspond
exactly to normal forms.  Figure \ref{fig:leafexp} is an example of
this correspondence, and more details can be found in
\cite{cfp,ctcomb,blake:diss}.

If an exposed caret has leaves numbered $i$ and $i+1$, then leaf $i+1$ must have leaf
exponent zero, since it is a right leaf.  If both trees in a tree pair diagram have exposed carets with
leaves numbered $i$ and $i+1$, then the corresponding normal form, computed via leaf exponents, contains
the generators $x_i$ to both positive and negative powers, but no instances of the generator $x_{i+1}$.
This is precisely the situation when the normal form can be reduced by a relator of $F$.  Thus the
condition that the normal form is unique is exactly the condition that the tree pair diagram is reduced.
This correspondence will be extended to elements of $T$ in the next section.

\subsection{Tree pair diagrams for elements of $T$}

We now discuss the relationship between words in $T$ and tree pair diagrams. This relationship is more
complicated in $T$ than it is in $F$.
%In $T$, it is no longer true that every tree pair diagram
%corresponds to an expression of the form
%$$
%x_{i_1}^{r_1}x_{i_2}^{r_2}\ldots x_{i_n}^{r_n}x_{j_m}^{-s_m}\ldots
%x_{j_2}^{-s_2}x_{j_1}^{-s_1},
%$$
%which is reduced if and only if the tree pair diagram is reduced.
%The relationship between tree pair diagrams and words in $T$ is
%not as pristine as it is in $F$.
%An element of $T$ is a piecewise linear homeomorphism of the
%circle, which can be represented by a pair of subdivisions of the
%interval $[0,1]$, as long as we indicate which interval of the
%target tree is the image of the first interval of the source tree,
%that is, the image of zero.  When these subdivisions are viewed as
%binary rooted trees in a tree pair diagram, we mark the image of
%$0$ on the target tree.
The representation of elements of $T$ by marked tree pair diagrams suggests a way to decompose an element
of $T$ into a product of three elements: the positive and negative parts together with a torsion part in
the middle, as described in \cite{cfp}.

\begin{defn} \label{factorization}
Let the marked tree pair diagram $(T_-,T_+)$ represent $g \in T$. If~~$T_-$ and $T_+$ each have $i+1$
carets, then we let
 $R$ be the all-right tree which has
 $i+1$ carets.  We can write $g$ as a product $pc_i^jq$, where:
\begin{enumerate}
 \item $p$, a positive word in the generators of $F$, is the normal form for the element of $F$ with tree pair diagram $(R,T_+)$, ignoring the marking on $T_+$.
 \item $c_i^j$ is a cyclic permutation of the leaves of $R$, with $1 \leq j \le i+2$, and
\item  $q$, a negative word, is the normal form for the element of $F$ represented by $(T_-,R)$.
\end{enumerate}
 Then the word $g=pc_i^jq$ is called the \emph{pcq factorization} of $g$ associated to the marked tree pair diagram $(T_-,T_+)$.
 In the special case
where $g \in F \subset T$, the $pcq$ factorization will just be the usual $pq$ factorization, as we
consider the $c$ part of the word to be empty (or equivalently, we can allow the exponent $j$ in the torsion part to be zero.)
\end{defn}

\begin{figure}
\includegraphics[width=5.2in]{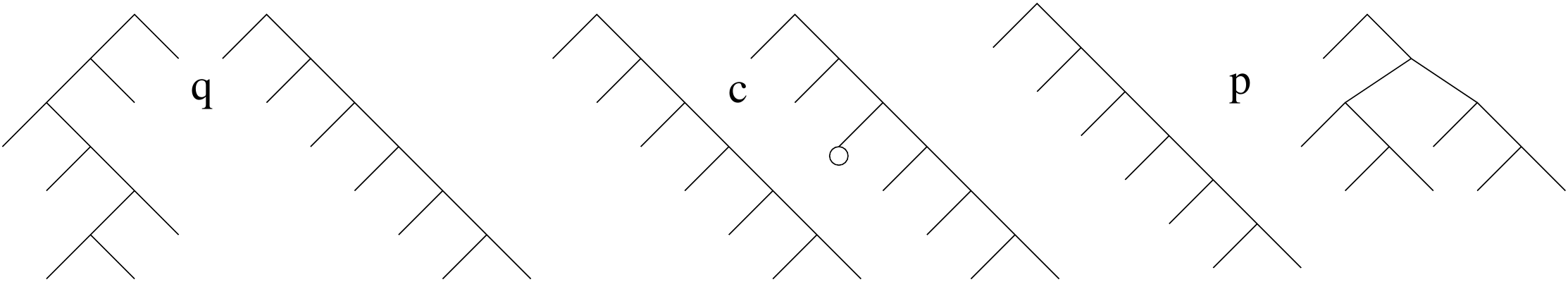}\\
\caption{ Three tree pair diagrams representing the word $x_1 x_2  c_5^5  x_2^{-2} x_1^{-1} x_0^{-2}$
factorized as $pcq$.\label{fig:threeparts}}
\end{figure}

Figure \ref{fig:threeparts} illustrates an example of an element of $T$ decomposed in this way.

The following theorem follows from the existence of these decompositions, and an algebraic proof of this
result is found in \cite{cfp}.

\begin{thm}[\cite{cfp}, Theorem 5.7]\label{pcq}Any element $x\in T  $ admits an expression of the form
$$
x_{i_1}^{r_1}x_{i_2}^{r_2}\ldots x_{i_n}^{r_n}\,c_i^j\,x_{j_m}^{-s_m}\ldots x_{j_2}^{-s_2}x_{j_1}^{-s_1},
$$ where $0 \leq i_1<i_2< \cdots <i_n$ and $0 \leq j_1<j_2<\cdots<j_m$ and either $1 \leq j<i+2$ or $c_i^j$ is not present.\end{thm}

We refer to any word satisfying the hypotheses of Theorem \ref{pcq} as a word in $pcq$ form for an
element of $T$ (just as words of this form with no $c_i^j$ term are called  words in $pq$ form in the
group $F$). Neither proof of the existence of $pcq$ forms gives an easy explicit method for transforming
a general word in the generators $x_i^{\pm 1},c_i$ into $pcq$ form without resorting to drawing tree pair
diagrams, so we will outline  an algebraic method below. We recall that the five types of relators we are using
in $T$ are:

\begin{enumerate}
\item $x_jx_i=x_ix_{j+1}$, if $i<j$ \item $x_kc_{n+1}=c_{n}x_{k+1}$,
if $k<n$ \item $c_nx_0=c_{n+1}^2$ \item $c_n=x_nc_{n+1}$  \item$c_n^{n+2}=1$
\end{enumerate}

%The proof of this result in \cite{cfp} is hard to implement, so we
%would like to give an explicit method for transforming a word in the %generators $x_i^{\pm 1},c_i$ into $pcq$ form.
%representing an element
%of $T$ in the $pcq$ form.
%For this, we need the following lemma.

\begin{lemma} [Pumping Lemma] The generators $x_i$ and $c_j$ of $T$
satisfy the following identities
$$
c_n^m=x_{n-m+1}c_{n+1}^m\qquad\qquad
 c_n^m=c_{n+1}^{m+1}x_{m-1}^{-1}
$$
if $1\le m < n+2$.
\end{lemma}

\begin{proof} This follows immediately from the relators. For instance, for the first identity, we have that
$$
c_n^m=c_n^{m-1}c_n=c_n^{m-1}x_nc_{n+1}
$$
by an application of relator of type (4). Now, several repeated applications
of relator (2) allow the $x_n$ to switch with the $c_n^{m-1}$ to
obtain the desired result. The second identity is the first one
taking inverses, and by noticing that $c_n$ has order $n+2$, we avoid
negative exponents for the $c$.
\end{proof}

We consider a word $w \in T$ written in the generators
$\{x_i,c_j\}$, and we describe explicitly an algebraic method of rewriting
it in $pcq$ form. The idea is to first combine occurrences of multiple $c_i$ generators into a power of a single one, and then to move the $x_n$ generators to the appropriate side of it. Consider first a
subword of the original word $w$ of type
$$
c_n^m\,w(x_i)\,c_k^l,
$$
where $w(x_i)$ is a word on the generators $x_n$ only, and which may
possibly be empty. We will apply relators to reduce this subword to
a word of the form $w_1(x_i) c_j^h w_2(x_i)$, where $w_1$ has only positive powers of $x_n$ generators and $w_2$ consists of only negative ones.
By the relators of type (1), we can assume that $w$ is of the form $pq$, that
is, with all positive powers of generators on the left and in increasing order of 
index, and all negative ones on the right and in decreasing order of index.
The goal is to move all the positive powers of $x_n$ generators to the left of
$c_n^m$ and all negative ones to the right of $c_k^l$. Although these moves may change the indices and powers of the $c_i$ generators, they merely change a power of a single $c_j$ generator to another power of a different single $c_k$ generator. To move all
the positive powers of generators to the left of $c_n^m$, we only need to use
relators of the type (2), assuming the index of $c$ is high enough.
If it is not, by repeated applications of the first identity of the
pumping lemma, the index can be increased arbitrarily, adding only
positive powers of generators to the {\it left\/} of $c_n^m$. When the subindex
is high enough, we can use relators of type  (2) to move all positive powers of generators of
$w$ past $c_n^m$. 
We note that a relator of type (3) may allow us to eliminate a occurance of $x_0$ to
the immediate right of $c_n^m$. It may be necessary
to combine the $c_i$ and $c_j$ generators obtained into a single term after this
elimination of $x_0$, as we see in an example:
$$
c_4^3x_1=c_4^2x_0c_5
$$
At this point $x_0$ cannot be moved farther, but we can use relator (3) to
obtain
$$
c_4c_5^3=x_5c_5^4
$$
with the last equality being an application of the pumping lemma to
$c_4$. We have achieved the goal of moving a positive power of a generator to the
left of $c_n^m$.

Moving the negative powers of the $x_n$ generators to the left is comparable. Using the second identity in the pumping lemma, we can
increase the index in $c_k^l$ as much as necessary to be able to
move all negative powers of the $x_n$ generators in $w$ to the right of $c_k^l$ using the relators (2) rewritten as $c_{n+1}x_{k+1}^{-1}=x_k^{-1}c_n$.
After this process, we will have a word consisting of positive powers of  $x_n$
generators, two powers of $c_i$ generators, and negative powers of $x_n$ generators. We now
 combine the powers of the two $c_i$ generators into a power of a single generator, by increasing
the smaller index to reach the larger. To do this, if the smaller
is on the left, we can use the first identity in the pumping lemma, and if
it is on the right, we can use the second one. This way no $x_n$ generator
will be added in between the two $c_i$ generators and after they have
the same index they can be combined into a power of a single generator. The positive powers of the $x_n$ generators now appear only to the left of the single power of the  $c_i$ generator, and negative
powers of  $x_n$ generators only to the right.

After repeated applications of this process to subwords of
the type $cwc$, we will have all occurrences of the $c_i$ generators
combined into a power of a single one. Our original word is now of the type
$$
w_1(x_i)\,c_n^m\,w_2(x_i),
$$
and $w_1$ and $w_2$ may again be assumed, after using relators (1), to be in $pq$ form. We only need to move the positive powers of generators in $w_2$ to the left of $c_n^m$ and the
negative powers of  $x_n$ generators of $w_1$ to the right of $c_n^m$, still maintaining a power of a single $c_i$ generator in the middle. We describe above as  the first step in our algorithm  precisely  how to do this. Furthermore, if the pumping lemma is needed to move a positive power of a generator to the left, recall that new positive positive powers of generators may appear in the word, but only to the left of the power of the $c_i$ generator.  
Hence, after moving each positive power of a generator,
all positive powers of generators in the word are to the left of $c_n^m$.

We now move each negative power of a generator to the right, and notice that the only cost of this is to add more negative powers of $x_n$ generators to the right of  $c_n^m$. When this is finished, the word has only positive powers of generators to the left of a power of a single $c$ and negative ones to the right.
Once the positive powers of the generators are together on the left side of the single $c$ term, 
we can reorder them if necessary using relators  of type (1), and similarly we can reorder the negative
part as well.

We will work an example as an illustration. Consider the word
$$
x_0^{-1}c_1x_3c_3^2x_1^{-1}
$$
The process starts by trying to move the $x_3$ to the left of
$c_1$. Since the index of $x_3$ exceeds the index of $c_1$, we cannot
apply a relator of type (2) directly. Using the pumping lemma, we write
$c_1=x_1c_2=x_1x_2c_3$. Hence our word is now the following, and
we can apply the relator of type to $c_3x_3$, obtaining:
$$
x_0^{-1}x_1x_2c_3x_3c_3^{2}x_1^{-1}=x_0^{-1}x_1x_2^2c_4c_3^2x_1^{-1}.
$$
We need to merge $c_4c_3^2$ into a single $c_i$ term. We increase the index of $c_3$ via
$c_3^2=c_4^3x_1^{-1}$ to obtain
$$
x_0^{-1}x_1x_2^2c_4c_4^3x_1^{-2}=x_0^{-1}x_1x_2^2c_4^4x_1^{-2}.
$$
The last step is to move the initial $x_0^{-1}$ to the
right side, using several relators of type (2)  to obtain
$x_0^{-1}c_4^4=c_5^4x_4^{-1}$. There is no need this time to
increase the index of $c_4^4$. The
final result is
$$
x_2x_3^2c_5^4x_4^{-1}x_1^{-2}
$$
which is in $pcq$ form.

The relationship between words in $pq$ form and tree pair diagrams in $F$ is different than the
relationship between $pcq$ forms and tree pair diagrams in $T$.  In  $F$, every tree pair diagram has a
$pq$ factorization associated to it, and any word in $pq$ form is in fact the $pq$ factorization
associated to a (not necessarily unique) tree pair diagram. Given any word in $F$ in $pq$ form, then we can form
a tree pair diagram for this element as follows. We consider reduced tree pair diagrams for $p$ and $q$,
and construct a tree pair diagram for the product $pq$ as described in Section \ref{multiplication}. The
middle trees of the four trees involved in the product are all-right trees. The all-right trees in this
decomposition may not have the same number of carets, so in forming the diagram for $pq$ we simply
enlarge the smaller of the two of these all-right trees (as well as the other tree in that diagram).
Since only right carets are ever added during this process, all of whose leaves have leaf exponent zero,
this results in a tree pair diagram whose $pq$ factorization is precisely the word $pq$ we began with.

In $T$, the correspondence between $pcq$ factorizations and general $pcq$ words is not as 
straightforward as in $F$. There is a difference between $pcq$ factorization and $pcq$ algebraic form.
Though  every element has a tree pair diagram  corresponding to a  $pcq$ factorization associated to it, there are words in algebraic $pcq$ form which are not the $pcq$ factorizations associated to a tree pair diagram.
The difficulty arises when the tree pair diagram for $c$ does not have as many carets as those for $p$ or
$q$, as adding right carets to enlarge $c$ appropriately necessitates adding generators to the normal
forms for $p$ and $q$, so the tree pair diagram one obtains by multiplying as in $F$ will not necessarily
have the original word as its factorization. For example, the word $x_1c_1$ is in algebraic  $pcq$ form,
yet it is not the $pcq$ factorization associated to any tree pair diagram. There is a different
representative for this element of $T$ which is the $pcq$ factorization associated to the reduced tree
pair diagram for this group element: $x_1c_2x_1^{-1}$. We prefer to work with words which are $pcq$
factorizations associated to tree pair diagrams, which will lead us to unique normal forms.

We can  algebraically characterize the words of type $pcq$ which are $pcq$ factorizations associated to
tree pair diagrams. The important condition is that the reduced tree pair diagram for $c$ should have at
least as many carets as those for $p$ and $q$.
%However, it may be
%possible to reduce the diagrams for $p$ and $q$, but not the
%diagram for $c$.
We say that words in $T$ with this property satisfy the {\em factorization condition}.

\begin{thm}
\label{thm:factor}  For elements in $T$ which are not in $F$, the word
$$
x_{i_1}^{r_1}x_{i_2}^{r_2}\ldots x_{i_n}^{r_n}\,c_i^j\,x_{j_m}^{-s_m}\ldots x_{j_2}^{-s_2}x_{j_1}^{-s_1}
,$$  where $0 \leq i_1<i_2< \cdots < i_n$, $0\leq j_1<j_2<\cdots<j_m$, and $1 \leq j < i+2$, is  the $pcq$
factorization associated to a tree pair diagram if and only if the number of carets in the reduced tree
pair diagram for $c_i^j$ is greater than or equal to the number of carets in the reduced tree pair
diagram for both of those for  the words $ x_{i_1}^{r_1}x_{i_2}^{r_2}\ldots x_{i_n}^{r_n}$ or
$x_{j_m}^{-s_m}\ldots x_{j_2}^{-s_2}x_{j_1}^{-s_1} $ in $F$.
\end{thm}

\begin{proof}
Given a tree pair diagram, by construction, the $pcq$ factorization associated to it satisfies the
factorization condition.  Given a word that satisfies the factorization condition, we can easily
construct the corresponding tree pair diagram as described above. The factorization condition
ensures that to perform the mulitiplication, $p \cdot c \cdot q$ as tree pair diagrams, it is only
necessary to (possibly)  add carets to the tree pair diagrams for the words $p$ and $q$.
This will not alter the normal form, and thus the 
 diagram constructed will indeed have the original word as its $pcq$
factorization.
\end{proof}

We can  compute  the number of carets of a reduced tree pair diagram for a  word $w \in F$ algebraically from the normal form of $w$,
as described by Burillo, Cleary and Stein in \cite{bcs}.

\begin{prop} [Proposition 2 of \cite{bcs}] \label{prop:ncarets}
Given a positive word in $w \in F$ in normal  form
$$
w=x_{i_1}^{r_1}x_{i_2}^{r_2}\ldots x_{i_n}^{r_n},
$$
then the number of carets $N(w)$ in either tree of a reduced tree diagram representing $w$ is
$$
N(w)=\max\{i_k+r_k+\ldots+r_n+1\} ,\text{ for }k=1,2,\ldots,n.
$$
\end{prop}

We can always  decide algebraically whether  $w \in T$, written in $pcq$ form,
corresponds to a tree pair diagram. We
 use Proposition \ref{prop:ncarets} to count the carets for the
positive and negative parts of the word.
The number of carets in a tree pair diagram for $c_i^j$ is
equal to $i+1$.

\section{Normal forms in $T$}

In $T$, we will declare the words in $pcq$ form which are $pcq$ factorizations associated to reduced
diagrams to be the normal forms for elements of $T$, similar to the approach used in $F$. However, it is
no longer true that these words cannot be shortened by applying a relator. As we saw with the normal
form $x_1c_2x_1^{-1}$ in $T$, a word may be the shortest word representing an element which satisfies the
factorization condition, yet there may be shorter words we can obtain by applying a relator which do not
satisfy the factorization condition.

%Since not all elements of $T$ arise from tree pair diagrams, the
%process of finding normal forms for elements of $T$ is slightly
%different than that for elements of $F$.  In $F$, each element has
%a unique normal form which corresponds to a unique reduced tree
%pair diagram.
%This difference between $F$ and $T$ appears when we study the
%effect of the relators on a given word.  In $F$, applying a
%relator to shorten a word corresponds to eliminating a superfluous
%caret from the tree pair diagram.  In $T$ the situation is more
%complicated.

Thus, when algebraically characterizing the normal form for elements of $T$, we restrict ourselves to
words of $pcq$ form which satisfy the factorization condition, regardless of whether or not a relator may
reduce the length of the word.  We next specify algebraic conditions which characterize the $pcq$
forms that correspond to normal forms, since we have given geometric conditions in Theorem \ref{thm:factor}

\begin{thm}\label{reductions} Let $w$ be a $pcq$ factorization for an element
$g \in T$ associated to a marked tree pair diagram in which each tree has  $i+1$ carets, where the $c$
part of the word  is $c_i^{j}$ with $1 \leq j < i+2$. A reduction of a pair of carets from the tree pair
diagram occurs only if the word $w$ satisfies one of the following conditions:
\begin{itemize}
\item[(1)] The pair of generators $x_{k-j}$ and
$x_k^{-1}$ appear, with $j \leq k< i$, and neither of the two generators
$x_{k-j+1}$ and $x_{k+1}^{-1}$ appear. The reduction corresponds to
applying the relator
$$
x_{k-j}c_{i}^jx_k^{-1}=c_{i-1}^j
$$ after applying relators from $F$ in the $p$ and $q$ parts of the word, if necessary, to make $x_{k+j}$ and $x_k^{-1}$ adjacent to $c_i^j$.
\item[(2)] The generator $x_{i-j}$ appears, and
$x_{i-j+1}$ does not. The reduction corresponds to applying
$$
x_{i-j} c_{i}^j=c_{i-1}^j
$$ after possibly using relators from $F$ as in (1).
\item[(3)] The pair of generators $x_{k+i-j+2}$ and
$x_k^{-1}$ for $0\leq k< j-2$ appear and neither one of the generators
$x_{k+i-j+1}$ or $x_{k+1}^{-1}$ appear. The reduction corresponds to
applying $$ x_{k+i-j+2}c_{i}^jx_k^{-1}=c_{i-1}^{j-1}
$$ after possibly applying relators from $F$.
\item[(4)] The generator $x_{j-2}^{-1}$ appears, and the generator
$x_{j-1}^{-1}$ does not appear. The reduction corresponds to
$$
c_{i}^{j} x_{j-2}^{-1}=c_{i-1}^{j-1}
$$ after possibly applying relators from $F$.
\end{itemize}
\end{thm}

\begin{figure}[b]
\includegraphics[width=5.82in]{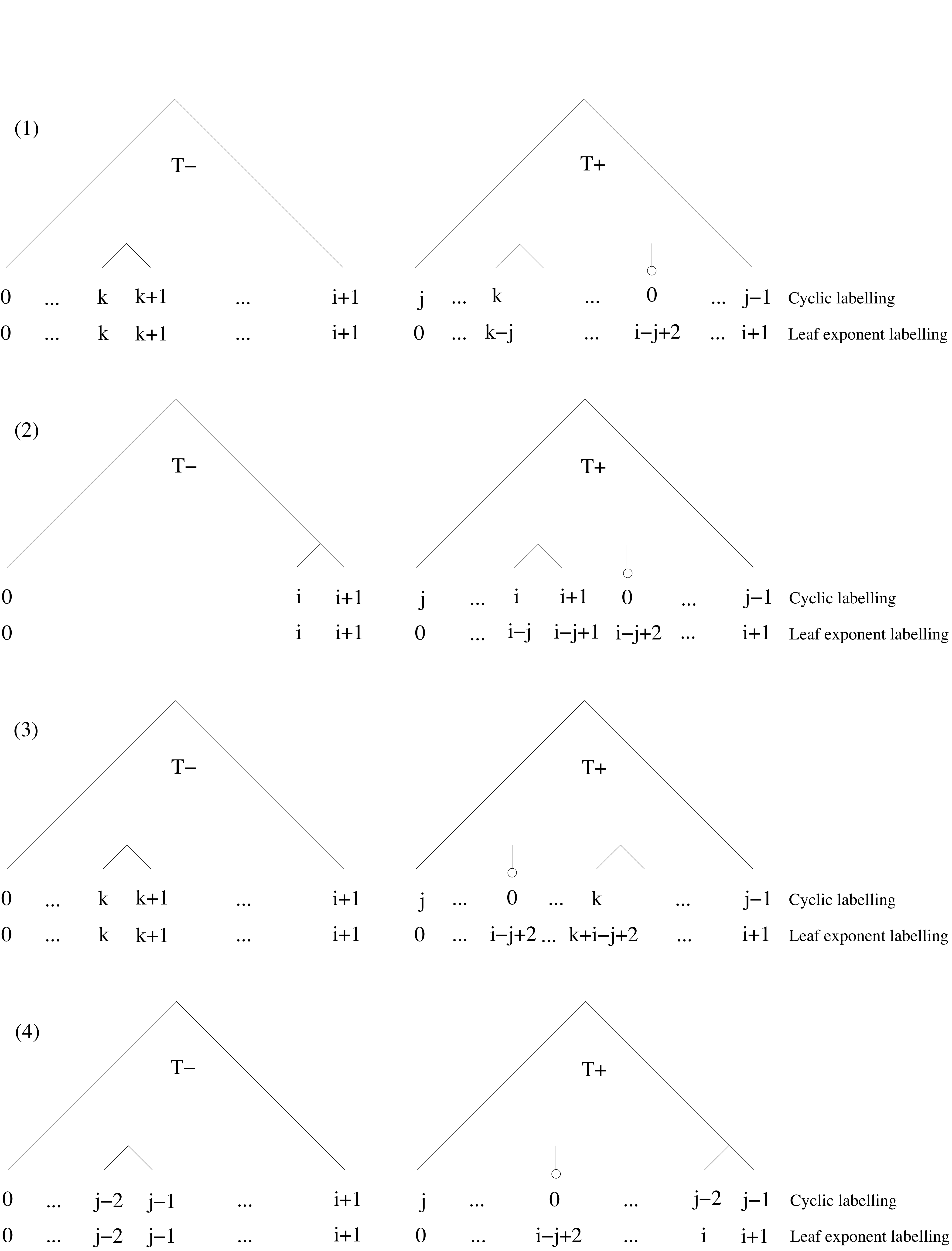}\\
\caption{The four cases in Theorem \ref{reductions}, showing the two
labellings on the leaves of the trees, the cyclic labelling which
indicates the correspondence of the leaves, and the leaf exponent
labelling which indicates the corresponding generators in the normal
form. \label{fig:normal}}
\end{figure}

\begin{proof}
Let $g \in T$ be represented by a marked tree pair diagram $(T_-,T_+)$.  If both trees have an exposed caret
whose leaves are identically numbered, then we call that a {\em reducible caret pair}, as it must be removed
in order to obtain the reduced tree pair diagram representing $g$.  We now consider algebraic conditions
corresponding to a reducible caret in a tree pair diagram.

In the tree pair diagram $(T_-,T_+)$ for $g \in T$,  there are two
ways of labelling the leaves in the target tree $T_+$.  The first
labelling corresponds to the order in which the intervals in the
subdivisions determined by these trees are paired in the
homeomorphism, and is called the cyclic labelling. The cyclic
labelling gives the marked leaf in the target tree the number zero,
and the other leaves are given increasing labels from left to right
around the leaves of the tree. The second labelling ignores the
marking and puts the leaves in increasing order from left to right,
beginning with zero.  The first labelling is used to determine which
leaves in $T_-$ are paired with which leaves in $T_+$, and the
second labelling is used in the computation of leaf exponents to
determine the powers of the generators that appear in the word.
 Figure \ref{fig:normal} shows the labellings  for
the four cases of the theorem.

Suppose that the tree pair diagram for $g \in T$ is not reduced.  The four cases above correspond to the
following four possible locations of a reducible caret relative to the marked leaf in the target tree.

\begin{itemize}
\item Case (1) of the thereom corresponds to the case when the left leaf of the reducible caret is
to the left of the marked leaf in $T_+$, but the reducible caret is
not the rightmost caret in $T_-$.
\item Case (2) corresponds to the special case when the reducible caret is a right caret in $T_-$, in which case necessarily its left leaf is to the left of the marked leaf in $T_+$.
Leaf exponents from leaves of right carets will always be zero and thus right carets cannot contribute generators
to the normal form.  They may still result in an exposed reducible caret, which occurs exactly in this
case, and the reduction will only affect
 the $q$ part of the
normal form. \item Case (3) corresponds to the case when the left leaf of the reducible caret is either
to the right of or coincides with the marked leaf in $T_+$, but the reducible caret is not the rightmost
caret in $T_+$. \item Case (4) corresponds to the special case when the reducible caret is a right caret
in $T_+$, in which case it cannot be to the left of the marked caret in $T_+$. As in Case (2), the
exposed caret in this case is a right caret and does not contribute a generator to the normal form, but
may still be reduced. This cancellation affects only the $p$ part of the normal form.
\end{itemize}

To see that these are all the possibilities, we note that $k$, the number of the left leaf in the cyclic numbering of the reducible caret in $T_-$, 
achieves all possible values in the cases above:
\begin{itemize}
\item If $0\le k<j-2$ we are in case (3).
\item If $k=j-2$ we are in case (4).
\item The case $k=j-1$ is impossible because the leaves $j-1$
and $j$ are at the two ends of the tree. With a cyclic ordering the
last and first leaves do not form a caret.
\item If $j\le k<i$ we are in case (1).
\item If $k=i$ we are in case (2).
\end{itemize}
Figure \ref{fig:normal} illustrates that these are
all the possibilities.\end{proof}

The conditions in Theorem \ref{reductions} together with the factorization condition algebraically
characterize our normal forms. The normal forms for elements in $F$ have already been characterized, so
we restrict to elements not in $F$ in our description.

\begin{thm} Any element $g\in T$ which is not an element of $F$
admits an expression of the form $pcq$ where
$$
p=x_{i_1}^{r_1}x_{i_2}^{r_2}\ldots x_{i_n}^{r_n}\qquad c=c_i^j\qquad q=x_{j_m}^{-s_m}\ldots
x_{j_2}^{-s_2}x_{j_1}^{-s_1},
$$  $0 \leq i_1<i_2< \cdots < i_n$, $0 \leq j_1<j_2<\cdots<j_m$, and $1 \leq j < i+2$.
Among all the words in this form representing an element, there is a unique one satisfying the following
conditions, which we call the normal form.
\begin{itemize}
\item The word satisfies the factorization condition, which we now state as $i+1\ge\max\{N(p),N(q)\}$.
\item The word does not admit any reductions, and thus its normal form
satisfies the following conditions:
\begin{itemize}
\item If there exists a pair of generators $x_{k-j}$ and
$x_k^{-1}$ simultaneously, for $j\le k< i$, then one of the
generators $x_{k-j+1}$ or $x_{k+1}^{-1}$ must appear as well.
\item If there is a generator $x_{i-j}$, then
$x_{i-j+1}$ must exist too.
\item If there exists a pair of
generators $x_{k+i-j+2}$ and $x_k^{-1}$ for $0\le k< j-2$, then one
of the generators $x_{k+i-j+1}$ or $x_{k+1}^{-1}$ must appear as
well.
\item If there exists a generator $x_{j-2}^{-1}$, then a generator
$x_{j-1}^{-1}$ must also appear.
\end{itemize}
\end{itemize}
\end{thm}

\begin{proof}

We claim that the conditions above precisely describe the set of unique normal forms for $T$. A $pcq$ word
satisfying the factorization condition is the $pcq$ factorization associated to a marked tree pair
diagram. However, if the $pcq$ word satisfies all four reduction conditions, we have just shown in the
previous theorem that this diagram is in fact the unique reduced diagram, and hence the word is in fact a
normal form. \end{proof}

We remark that the Pumping Lemma together with the reductions in Theorem \ref{reductions} give an
explicit way of algebraically transforming any word in the generators of $T$ into a normal form. 
Given any word,  we rewrite it in $pcq$ form using the process described following the Pumping Lemma. If
the resulting word does not satisfy the factorization condition, then we iterate the Pumping Lemma until
we obtain a word for which the factorization condition is satisfied.  The Pumping Lemma increases the
number of carets for $c$ and the number of carets for one of the words $p$ and $q$.  Once a word is
obtained which satisfies the factorization condition, there must be a corresponding tree pair diagram for
the element. Now, if the word satisfies any of the reduction conditions in Theorem \ref{reductions}, we
apply them successively using the relators described there. This method thus produces the unique normal
form.

\section{The word metric in $T$}
\subsection{Estimating the word metric}
For metric questions concerning $T$, we must consider a finite generating set instead of the one used to
obtain the normal form for elements.  We now approximate the word length of an element of $T$ with
respect to the generating set $\{x_0,x_1,c_1\}$, using information contained in the normal form and the
reduced tree pair diagram.  These estimates are similar to those for the estimates of word metric in $F$ with respect to the
generating set $\{x_0,x_1\}$ described  \cite{burillo}, \cite{bcs}.

\begin{thm} \label{thm:D} Let $w \in T$ have normal form
$$
w=x_{i_1}^{r_1}x_{i_2}^{r_2}\ldots x_{i_n}^{r_n}\,c_i^j\,x_{j_m}^{-s_m}\ldots
x_{j_2}^{-s_2}x_{j_1}^{-s_1}.
$$
We define
$$
D(w)=\sum_{k=1}^nr_k+\sum_{l=1}^ms_l+i_n+j_m+i.
$$
Let $|w|$ denote the word metric in $T$ with respect to the generating set $\{x_0,x_1,c_1\}$.  There
exists a constant $C>0$ so that for every $w \in T$,
$$
\frac{D(w)}C\le|w|\le C\,D(w)
$$
and similarly, for $N(w)$ the number of carets in the reduced tree pair diagram representing $w$,
$$
\frac{N(w)}C\le|w|\le C\,N(w).
$$
\end{thm}

\begin{proof}
These inequalities follow from the correspondence between the normal form and the tree pair diagram for
an element $w \in T$.  It is clear, from Proposition \ref{prop:ncarets}, that $ N(w)\ge \sum_{k=1}^nr_k$,
$N(w)\ge \sum_{l=1}^ms_l$, $N(x)\ge i_n$, and $N(w) \ge j_m$. The inequality $ N(w)\ge i $ is clear from
the fact that $c_i$ has $i+1$ carets. These inequalities prove that
$$
D(w)\le 5\,N(w).
$$

We rewrite the generators $x_i$ and $c_j$ in terms of $x_0$, $x_1$  and $c_1$ and look at the lengths of
the resulting words to obtain the inequality
$$
|w|\le C\,D(w)
$$
for some constant $C>0$.  Combining the two inequalities above, we have
$$
|w|\le C'\,N(w).
$$

To obtain lower bound on the word length, we consider the fact that the tree pair diagram for each
generator has either two or three carets. If $u$ is a word in $x_0$, $x_1$ and $c$ with length $n$, then
as these generators are multiplied together, each product may add at most $3$ carets to the tree pair
diagram. Thus the diagram for $u$ will have at most $3n$ carets.  It then follows that
$$
N(w)\le 3|w|.
$$
Combining this with the above inequality, we obtain the desired bounds.
\end{proof}

We use Theorem \ref{thm:D} to show that the inclusion of $F$ in $T$ is a quasi-isometric embedding.  This
means that there are constants $K>0$ and $C$ so that for any $w,z \in F$ we have
$$\frac{1}{K} d_{F}(w,z) - C \leq d_T(w,z) \leq Kd_F(w,z) + C$$
where $d_F$ and $d_T$ represent the word metric in $F$ and $T$ respectively, with regard to the
generating set $\{x_0,x_1\}$ of $F$ and $\{x_0,x_1,c_1\}$ of $T$.

When considering whether the inclusion of a finitely generated subgroup $H$ into a finitely generated
group $G$ is a quasi-isometric embedding, we can instead  equivalently show that the distortion function
is bounded.  The distortion function is defined by
$$h(r) = \frac{1}{r} \max \{|x|_H :  x \in H \mbox{~and~}  |x|_G \leq r\}.$$

Word length in $F$ is comparable to the number of carets in the reduced tree pair diagram representing
the word, by Theorem 3 of \cite{bcs} or more directly by Fordham's method \cite{blakegd}.  This, combined with Theorem \ref{thm:D} easily shows that the
distortion function is bounded, and thus proves the following corollary with respect to
the generating sets $\{x_0, x_1\}$ and $\{x_0,x_1,c_1\}$ and thus all pairs of finite generating sets:

\begin{cor} \label{cor:qiemb}
The inclusion of $F$ in $T$ is a quasi-isometric embedding.
\end{cor}

\subsection{Comparing word length in $F$ and $T$\label {isomembed}}

Although Corollary \ref{cor:qiemb} shows that $F$ is quasi-isometrically embedded in $T$, in fact the
word length of many elements of $F$ does not change at all when these elements are considered as elements of
$T$, with respect to natural finite generating sets. As an example of this phenomenon, we characterize
one type of element of $F$ whose word length is unchanged when viewed as an element of $T$, using the
generating set $\{x_0,x_1\}$ for $F$ and $\{x_0,x_1,c_0\}$ for $T$.  These are elements $w \in F$ for
which $N(w)$ exceeds the word length $|w|_F$.  Fordham \cite{blakegd} computes $|w|_F$ by assigning an
integer weight between zero and four to each pair of carets in the tree pair diagram representing $w$. In
a given word there are at most two weights of zero.  Here we investigate words in which most weights are
one.  Such words, for example, are represented by tree pair diagrams with no interior carets having right
children.

\begin{thm} \label{thm:wordlength}
If $w \in F$ with $N(w) \geq |w|_F + 1$ then $|w|_T = |w|_F$, where word length if computed with respect
to the generating set $\{x_0,x_1\}$ for $F$ and $\{x_0,x_1,c_0\}$ for $T$.
\end{thm}

This theorem is proved by taking a word in the generators of $T$, and analyzing how each generator
changes the intermediate tree pair diagram as one builds up the final tree pair diagram for $w$.
Carefully controlling the process allows one to obtain an upper bound on $N(w)$ in terms of the
length of the word. If the word is actually shorter than $|w|_F$, then this bound, considered together
with the lower bound given by the hypothesis, yields a contradiction. We immediately obtain the
following corollary, since $|x_0^n|_F = |x_1^n|_F = n$, while $N(x_0^n) = n+1$ and $N(x_1^n) = n+3$.

\begin{corollary}
The elements $x_0^n$ and $x_1^n$ have word length $n$ in both $F$ and $T$ with respect to the finite
generating sets $\{x_0,x_1\}$ and $\{x_0,x_1,c_0\}$ respectively.
\end{corollary}

\section{Torsion elements}
Although the group $F$ is torsion free, both $T$ and $V$ contain torsion elements. It is easy to
construct torsion elements in $T$ or $V$ by choosing any binary tree $S$ and making any marked tree pair
diagram with $S$ as both source and target tree.  If the labelling of the target tree is the same as the
labelling of the source tree, we get an unreduced representative of the identity; otherwise, we get a
non-trivial torsion element. If this is an element of $T$, the tree pair diagram has $pcq$ factorization
in which $q=p^{-1}$.
   In fact, any torsion element can be represented by such a tree pair diagram, though its
 reduced marked tree pair diagram may well not have the same source and target trees, corresponding to the fact
 that although it has a $pcq$ word where $q=p^{-1}$, the normal form may well not have this special balanced appearance.

\begin{prop}\label{p:torsion}
If $f\in F,T$ or $V$  is a torsion element, then  it can be represented by a (marked) tree pair diagram
with the same source and target trees.
\end{prop}
Before proving Proposition \ref{p:torsion}, we establish some notation which links the analytic and
algebraic representations of these groups.  For $f \in F$, $T$, or $V$, if $(T_-,T_+)$ is a marked
tree pair diagram representing $f$, then it is sometimes convenient to denote the tree $T_+$ by $f(T_-)$.
The tree $T_-$ corresponds to a certain subdivision of the circle, which maps under $f$ linearly to
another subdivision of the circle. This subdivision is represented by the tree $T_+$, and the marking
describes where each subinterval of the circle is mapped. The element $f$ can be thought of as mapping
the leaves of $T_-$ to the leaves of $f(T_-)=T_+$, where the marking defines this mapping of the leaves.
If $f$ does not have a tree pair diagram in which the tree $T$ appears as the source tree, then the
symbol $f(T)$ has no meaning.

Given two rooted binary trees $T$ and $T'$, we say that $T'$ is an \emph{expansion} of $T$ if $T'$ can be
obtained from $T$ by attaching the roots of additional trees to some subset of the leaves of $T$. We
observe that if $(T, f(T))$ is a marked tree pair diagram for $f$, and $T'$ is an expansion of $T$, then
there is always a tree pair diagram $(T', f(T'))$ for $f$, and $f(T')$ is an expansion of $f(T)$. Given
two rooted binary trees $S$ and $T$, by the \emph{minimal common expansion} of $S$ and $T$ we mean the
smallest rooted binary tee which is an expansion of both $S$ and $T$.  Using this language, if $(T,f(T))$
and $(S,g(S))$ are marked tree pair diagrams for $f$ and $g$ respectively, the process described in
Section 2.3 for creating a tree pair diagram for the product $gf$ could be summarized as follows. If $E$
is the minimal common expansion of $f(T)$ and $S$, then there are tree pair diagrams $(f^{-1}(E),E)$ for
$f$, $(E,g(E))$ for $g$, and $(f^{-1}(E),g(E))$ for $gf$ (with appropriate markings).

\begin{proof}
Suppose that $f$ is a torsion element.  We begin by describing the construction of (marked) tree pair
diagrams $(A_n,B_n)$ for $f^n$ for every $n \geq 1$.  These tree pair diagrams are constructed
inductively, viewing $f^n$ as a product $(f^{n-1})(f)$. For $n=1$, let $(A_1,B_1)$ be the reduced marked
tree pair diagram for $f$. Throughout this procedure, although markings are carefully carried through in
either $T$ or $V$, since our goal is merely to produce a tree pair diagram for $f$ with the same source
and target trees (regardless of marking), only the trees themselves are relevant for this argument. Hence
we suppress mention of any markings throughout the construction. If $k \geq 2$, suppose the marked tree
pair diagram $(A_{k-1},B_{k-1})$ for $f^{k-1}$ has been constructed. Let $E_{k-1}$ be the minimal common
expansion of the trees $A_1$ and $B_{k-1}$. Then $f^k$ has tree pair diagram
$(f^{-(k-1)}(E_{k-1}),f(E_{k-1}))$, and we let  $B_k=f(E_{k-1})$ and $A_k=f^{-(k-1)}(E_{k-1})$.

By construction, $A_{k+1}$ is an expansion of $A_k$ for all $k \geq 1$. We claim also that $B_{k+1}$ is
an expansion of $B_k$ for all $k \geq 1$. For $k=1$, $E_1$ is by definition an expansion of $A_1$, which
implies that $B_2=f(E_1)$ is an expansion of $B_1=f(A_1)$. Suppose inductively that $B_k$ is an expansion
of $B_{k-1}$. Now $E_k$ is an expansion of $B_k$ and $A_1$, so $E_k$ is an expansion of $B_{k-1}$ and
$A_1$. But $E_{k-1}$ is the minimal common expansion of $B_{k-1}$ and $A_1$, so $E_k$ is an expansion of
$E_{k-1}$, which implies that $B_{k+1}=f(E_k)$ is an expansion of $B_k=f(E_{k-1})$.

Since there exists a positive integer $m$ such that $f^m$ is the identity, it follows that all tree pair
diagrams for $f^m$ must have the same source and target trees. Hence $A_m=B_m$, and then since $A_m$ is
an expansion of $A_1$, $B_m$ is an expansion of $A_1$. But since $E_{m-1}$ is the minimal common
expansion of $B_{m-1}$ and $A_1$, the fact that $B_m$ is an expansion of both $B_{m-1}$ and $A_1$ implies
that $B_m=f(E_{m-1})$ is an expansion of $E_{m-1}$. But they have the same number of carets, so in fact
$f(E_{m-1})=E_{m-1}$. In other words, the tree pair diagram $(E_{m-1}, B_n=f(E_{m-1}))$ is the desired
tree pair diagram for $f$.
\end{proof}

\begin{cor} An element of $T$ is torsion if and only if it is a conjugate of some $c_i^j$.
\end{cor}

\begin{proof} If an element is torsion, then it admits a diagram with two equal trees. The $pcq$ factorization
associated with this diagram has the form  $pc_i^jp^{-1}$, where $p$ is a positive element of $F$.
\end{proof}

A particularly natural torsion subgroup is the subgroup $R$ of pure
rotations, where by a pure rotation we mean a rotation by
$d=\frac{a}{2^n}$ (where $a$ is not divisible by 2). Such pure
rotations were used in Section \ref{sec:rotationnum} to conjugate
the fixed point of a homeomorphism to 0.

This subgroup is isomorphic to the group of dyadic rational numbers modulo 1, which has a 2-adic metric
as follows: if $x=\frac{p}{2^l}$, $y= \frac{q}{2^m}$, and $z= |x-y|  = \frac{r}{2^k}$, where $p, q$ and
$r$ are odd, then $d(x,y)=2^k$. With respect to this metric, the subgroup of rotations is
quasi-isometrically embedded in $T$.

\begin{prop} The subgroup $R$ of the pure rotations, with the 2-adic metric, is quasi-isometrically
embedded in $T$.
\end{prop}

\begin{proof}
We note that if $g\in T$ is the rotation by $\frac{a}{2^n}$ where $a$ is not divisible by $2$, then there
are $2^{n}-1$ carets in the reduced tree pair diagram representing $g$, so $N(g)=2^{n}-1$. Since we have
shown that the word length of $g$ in $T$ is bi-Lipschitz equivalent to $N(g)$, the proposition follows.
\end{proof}

\bibliographystyle{plain}

\end{document}